\documentclass[12pt]{amsart}
\usepackage{amsfonts}
\newcommand{\numb}[1]{(\ref{#1})}
\newcommand{\sgn}{\mathop{\rm sgn}}

\newtheorem{thm}{Theorem}

\newtheorem{lemma}{Lemma}
\newcommand{\R}{{\mathbb R}}

\newcommand{\C}{{\mathbb C}}

\newcommand{\su}{{\mathfrak s}{\mathfrak u}}
\newcommand{\so}{{\mathfrak s}{\mathfrak o}}
\newcommand{\prodin}{\prod_{i=1}^n}
\newcommand{\prodkn}{\prod_{k=1}^n}
\newcommand{\sumweyl}{\sum_{w\in W}}
\newcommand{\sumhypoct}{\sum_{\sigma\in S_n}\sum_{\epsilon_i=\pm1}}
\newcommand{\btel}{b_t(\eta,\lambda)}
\newcommand{\btelp}{b_t(\eta,\lambda')}
\newcommand{\bptel}{b'_t(\eta,\lambda)}

\def\1/2{{1\over 2}} % One-half
\title[Brownian Motion in a Weyl Chamber]{Brownian Motion in a Weyl
Chamber,\\ Non-colliding Particles, and Random Matrices}
\author{David J. Grabiner}\thanks
{Supported by an NSF Graduate Fellowship and NSF Postdoctoral Fellowship,
and a postdoctoral fellowship from the Center for Discrete Mathematics
and Theoretical Computer Science at the Hebrew University of Jerusalem.
Some of this work was done while visiting the Mathematical Sciences
Research Institute; MSRI is supported by NSF Grant DMS-9022140.}

\email{grabiner@math.lsa.umich.edu}

\address{Department of Mathematics, University of Michigan, Ann Arbor,
MI 48109-1109}

\begin{document}

\begin{abstract}  Let $n$ particles move in standard Brownian motion in
one dimension, with the process terminating if two particles collide.
This is a specific case of Brownian motion constrained to stay inside a
Weyl chamber; the Weyl group for this chamber is $A_{n-1}$, the
symmetric group.  For any starting positions, we compute a determinant
formula for the density function for the particles to be at specified
positions at time $t$ without having collided by time $t$.  We show that
the probability that there will be no collision up to time $t$ is
asymptotic to a constant multiple of $t^{-n(n-1)/4}$ as $t$ goes to
infinity, and compute the constant as a polynomial of the starting
positions.  We have analogous results for the other classical Weyl
groups; for example, the hyperoctahedral group $B_n$ gives a model of
$n$ independent particles with a wall at $x=0$.

We can define Brownian motion on a Lie algebra, viewing it as a vector
space; the eigenvalues of a point in the Lie algebra correspond to a
point in the Weyl chamber, giving a Brownian motion conditioned never to
exit the chamber.  If there are $m$ roots in $n$ dimensions, this shows
that the radial part of the conditioned process is the same as the
$n+2m$-dimensional Bessel process.  The conditioned process also gives
physical models, generalizing Dyson's model for $A_{n-1}$ corresponding
to ${\mathfrak s}{\mathfrak u}_n$ of $n$ particles moving in a diffusion
with a repelling force between two particles proportional to the inverse
of the distance between them.
\end{abstract}
\maketitle

\section{Introduction}

Let $n$ particles move in standard Brownian motion in one dimension,
with the process terminating if two particles collide.  Given the
starting positions, we can use a reflection argument to calculate the
density function for the particles to be at specified positions at time
$t$ without having collided by time $t$.  Using this density function
and the theory of Lie algebras, we can prove the following results;
Theorems~\ref{processthm} and~\ref{forcethm} were first proved by
Dyson~\cite{dyson}.

\begin{thm}\label{nocollthm}
For any starting positions, the probability that there will be no
collision up to time $t$ is asymptotic to a constant multiple of
$t^{-n(n-1)/4}$ as $t$ goes to infinity; the constant is a known
polynomial in the starting positions.
\end{thm}

\begin{thm}\label{radialthm}
Given that there is no collision up to time $t$, the distribution of the
radius of the vector whose coordinates are the positions of the
particles, divided by the square root of $t$, converges in measure to
the distribution of the Bessel process with parameter $n(n+1)/2$ at time
$1$, which is the radial part of an $n(n+1)/2$-dimensional standard normal.
\end{thm}

\begin{thm}\label{processthm}
We can construct an $n$-dimensional Brownian motion which is conditioned
for no two particles ever to collide.  If we take Brownian motion on the
space of Hermitian matrices, the induced process on the eigenvalues is
the same process.  If the starting point is appropriately chosen at the
starting radius, the process given by the radial part of the conditioned
Brownian motion is the same process as the $n^2$-dimensional Bessel
process.
\end{thm}

\begin{thm}\label{forcethm}
This conditioned process is identical to the process obtained by $n$
particles moving in a one-dimensional diffusion with constant
infinitesimal variance, with a repelling force between two particles
proportional to the inverse of the distance between them; that is, its
infinitesimal generator has $\sigma^2_{ij}(\vec x,t)=\delta_{ij}$,
$\mu_i(\vec x,t)=\sum_{i\ne j}1/(x_i-x_j)$.
\end{thm}

This is a specific case of Brownian motion in a Weyl chamber; the vector
whose coordinates are the locations of the $n$ particles is constrained
to stay inside the chamber.  The Weyl groups is $A_{n-1}$, the symmetric
group.

We have similar results for other Weyl groups; Theorems
\ref{nocollthm} and~\ref{radialthm} have analogues for the classical
Weyl groups, and Theorems \ref{processthm} and~\ref{forcethm} have
analogues for all Weyl groups.  For example, the Weyl group $B_n$ models
$n$ particles in independent Brownian motion with an absorbing wall at
$x=0$.  In this case, the asymptotic probability is $ct^{-n^2/2}$ (and
the constant is again known), with the radial part corresponding to the
$n(n+1)$-dimensional Bessel process; the conditioned process has a
radial part which is the $n(2n+1)$-dimensional Bessel process.  This
generalizes the results of Pitman and Williams~\cite{pitman,williams}
that one-dimensional Brownian motion conditioned to stay positive is the
same process as the three-dimensional Bessel process.

Our reflection argument is a generalization of the {\em reflection
principle}, a standard argument in the analysis of both discrete random
walks and Brownian motion.  In the discrete case, it is used in the
classical formula for the Catalan numbers, which enumerate the
arrangements of $n$ $+1$'s and $n$ $-1$'s so that none of the partial
sums are negative.  Similarly, it can be used to study Brownian motion
in one dimension with an absorbing barrier at $x=0$ and a known starting
point~\cite{harrison}.  It has also been extended to multiple
reflections in one dimension to study Brownian motion with two absorbing
barriers~\cite{freedman}.

The reflection principle has been generalized to multiple dimensions.
For example, the {\em ballot problem}, a classical problem in random
walks, asks how many ways there are to walk from the origin to a point
$(\lambda_1,\ldots,\lambda_n)$, taking $k$ unit-length steps in the
positive coordinate directions while staying in the region $x_1\ge
x_2\ge\cdots\ge x_n$.  The solution is known in terms of the hook-length
formula for Young tableaux; a combinatorial proof, using a reflection
argument, is given in~\cite{WM,Z}.

The same reflection argument has also been applied to the case of $n$
independent diffusions, or discrete processes which cannot pass each
other without first colliding.  Using this method, Karlin and
McGregor~\cite{KM} give a determinant formula for the probability or
measure for the $n$ particles, starting at known positions, not to have
collided up to time $t$ and to be in given positions.  Hobson and
Werner~\cite{HW} generalize this argument to $n$ particles in an
interval or circle, and use this to prove a result analogous to
Theorem~\ref{forcethm} for $n$ particles on the circle.

Gessel and Zeilberger~\cite{GZ}, and independently Biane~\cite{biane92}
give a further generalization.  For certain ``reflectable'' random
walk-types, we can count the number of $k$-step walks between two points
of a lattice, staying within a chamber of a Weyl group, in terms of
numbers of unconstrained walks.  The steps must have certain allowable
lengths and directions.  In~\cite{decomp}, all cases in which this
method applies are enumerated, and determinant formulas are given for
many important cases, including walks in the classical Weyl chambers.

The argument of~\cite{GZ} can be generalized to Brownian motion in any
Weyl chamber or chamber of a Coxeter group, with either absorbing or
reflecting boundary conditions.  We prove this generalization, and then
use the result to compute determinant formulas for Brownian motion in
the Weyl chambers of $A_{n-1}$, $B_n=C_n$, and $D_n$.  The $A_{n-1}$ and
$B_n$ cases are applicable to the independent motion of $n$ particles in
one dimension.  The $A_{n-1}$ formula appears in~\cite{KM}, using the
model of $n$ independent particles rather than motion in a Weyl chamber.
The cases of the affine Weyl groups $\tilde A_{n-1}$ and $\tilde B_n$
are studied in~\cite{HW}, also as models of $n$ independent particles,
in a circle or interval.

These determinants factor into forms which can be easily analyzed; this
allows us to find the asymptotic probability that there will be no
collision up to time $t$ for these three cases, for any starting point.
In the case of $A_{n-1}$, we get a simple formula for the actual measure
as well.

%These asymptotics suggest a general technique for constructing Brownian
%motion conditioned to stay in a chamber of any finite Weyl or Coxeter
%group, using a harmonic transformation of standard Brownian motion.  We
%construct this motion, and the analogous discrete random walks.  The
%possibility of transforming Brownian motion in this way is mentioned
%in~\cite{biane}, referring to~\cite{BB} for a classification of the
%harmonic functions involved.

Weyl chambers arise naturally as the spaces of eigenvalues of elements
of Lie algebras.  Also, Brownian motion on a Lie algebra is symmetric
under an orthonormal change of basis.  There is thus a natural
correspondence between standard Brownian motion on the Lie algebra and
some diffusion on the Weyl chamber.  Dyson~\cite{dyson} showed that this
diffusion has the same generator as the conditioned Brownian motion for
$A_{n-1}$; this generalizes naturally to other Lie algebras.  We can
thus use the known properties of random matrices to study the
distribution.

DeBlassie~\cite{deblassie} uses a different approach to give a more
general formula for asymptotics and density functions for a general
class of cones, which include the cases discussed here.  A cone is
defined as the union of all rays from the origin which intersect the
unit sphere in a connected open set $C$.  The Laplace-Beltrami operator
$L_{S^{n-1}}$ on the unit sphere is the non-radial portion of the
ordinary Laplacian.  If $\lambda_1$ is the eigenvalue of $L_{S^{n-1}}$
with smallest absolute value on the space of all $L^2$ functions on $C$
which vanish continuously on the boundary of $C$, then the asymptotic
probability is a constant multiple of
\[
t^{-\left[-(n-2)+\sqrt{(n-2)^2+4\lambda_1^2}\right]/4}.
\]
This result gives asymptotics for a large class of cones.  The
coefficients depend on the eigenvalues of the Laplace-Beltrami operator
rather than on explicit coefficients; these eigenvalues are known for
Weyl chambers~\cite{BB}.  The asymptotics can thus be computed from
these formulas; in theory, the explicit density functions can also be
computed, but as infinite series, This formula also shows that the
asymptotics for all such cones exist and are powers of $t$, with no
other terms such as logarithms.

The asymptotic probability is known to be a constant multiple of
$t^{-m/2}$ for a wedge of angle $\pi/m$ in two
dimensions~\cite{burkholder}.  The result holds in general, although the
region is only a chamber of a Coxeter group (the dihedral group) if $m$
is an integer.

O'Connell and Unwin~\cite{OU} study the opposite
asymptotic problem to ours, computing an asymptotic for the probability
that $n$ particles in independent Brownian motion will have a collision
up to time $t$ when $t$ is small compared to the initial separation of
the particles.

Our results are organized as follows.  Section 2 contains the basic
definitions.  In section 3, we prove the basic reflection result, and in
section 4, we apply this result to get the determinant formulas.  In
section 5, we prove Theorems \ref{nocollthm} and~\ref{radialthm}, and
their analogues for $B_n$ and $D_n$.  In section 6, we construct the
conditioned motion, both by Lie theory and by $h$-transformation,
proving Theorems~\ref{processthm} and \ref{forcethm} and their
analogues, and use the $h$-transformation to find a physical model for
all finite Coxeter groups.

\section{Definitions}

We will study a process with continuous sample paths in $\R^n$, either
unconstrained or constrained by a chamber.  In the constrained case, we
may have either an absorbing boundary condition, causing the process to
terminate when it hits a wall, or a reflecting boundary condition.  All
references to the analogous discrete problem are discussed
in~\cite{decomp}.

We require that our chamber $C$ be a chamber of a finite or affine {\em
Coxeter group}.  In the finite case, $C$ is defined by a system of
simple roots $\Delta \subset \R^n$ as
\begin{equation}
 C = \{ \vec{x} \in \R^n \mid (\alpha, \vec{x})\geq 0 \mbox{ for all } 
\alpha \in \Delta \},
\end{equation}
and the orthogonal reflections $r_{\alpha}$: $\vec{x} \mapsto \vec{x} -
\frac{ 2(\alpha, \vec{x})}{(\alpha,\alpha)} \alpha$, generate a finite
group $W$ of linear transformations, the {\em Coxeter group}.  In the
affine case, the hyperplanes of reflection which define $C$ do not all
pass through the origin, and the group $W$ is infinite, but if $T$ is
the subgroup of all translations, $W/T$ must be finite.  In the
analogous discrete problem, the steps of the the random walk must
generate a lattice $L$ which is stable under the action of $W$; in this
case, $C$ is a {\em Weyl chamber} and $W$ a {\em Weyl group}.  

Let $X(t)$ be a Markov process with continuous sample paths with values
in $\R^n$; that is, the distribution of $X(t_2)$ given $X(t_1)$ is
independent of $X(t)$ at any $t$ outside the interval $[t_1,t_2]$.  We say
that the constrained motion is {\em reflectable} if the increments of
the unconstrained motion are symmetric under the Coxeter group; that is,
the distribution of $X(t_2)$ given $X(t_1)$ is the same as the
distribution of $w(X(t_2))$ given the $w$-image $w(X(t_1))$ as the
starting point for any $w\in W$.

Standard Brownian motion is reflectable for any Coxeter group.  For a
finite Coxeter group, in which all planes of reflection pass through the
origin, any diffusion with variance dependent only on time and the
radius, and drift dependent on time and symmetric with respect to
rotations and reflections about the origin, is reflectable; for example,
there could be an absorbing or reflecting barrier at $|\vec x|=R$.

As another example, consider the case in which each coordinate $x_i(t)$
is an independent identical diffusion; this could model $n$ independent
particles instead of one particle in $n$ dimensions.  If our Coxeter
group is the symmetric group $A_{n-1}=S_n$ (giving the chamber
$x_1>x_2>\cdots>x_n$), it permutes the particles, so the process is
reflectable under this action; this case is discussed in~\cite{KM}.  If
the individual diffusions are symmetric about $x_i=0$, then the product
process will also be reflectable under the hyperoctahedral group $B_n$,
which includes all permutations with any number of sign changes.

In the discrete case, reflectability requires the additional condition
that the walk cannot go from inside the chamber to outside it without
stopping on a wall~\cite{GZ}; this is our condition of continuous sample
paths, which is satisfied by any diffusion.

For fixed $t$, this process defines a probability measure
$P_t(A)=P\{X(t)\in A\}$ which represents the chance that this process,
if started at 0, will be in a set $A$ at time $t$.  We let $c_t(\vec x)$
be the density function of this probability measure with respect to
Lebesgue measure on $\R^n$; that is, 
\[
\int_{\vec x\in A} c_t(\vec x) = P_t(A).
\]

Now, we study the case in which the motion is constrained by a 
chamber, with either absorbing or reflecting boundary conditions.  We
must now fix the starting point $\eta$, since the process now depends on
it in a non-trivial way.  This generates new stochastic processes
$Y(\eta,t)$ for absorbing boundary conditions and $Y'(\eta,t)$ for
reflecting boundary conditions.  These give probability measures
$Q_t(\eta,A)$ and $Q'_t(\eta,A)$ which give the probability that the
process, if started at $\eta$, will be in the set $A$ at time $t$, and
these have density functions $\btel$ and $\bptel$.  Note that the total
measure $Q_t(\eta,\R^n)$ will be less than 1 because the process
terminates when it reaches a wall of the chamber.  The total measure
$Q'_t(\eta,\R^n)$ for reflecting boundary conditions will still be 1.

If our process is a Brownian motion, it has drift $\mu_i$ and variance
$\sigma_i^2$ in each coordinate direction.  If the Coxeter group
contains any reflection in the $x_i,x_j$-plane other than a sign change
of $x_i$ or $x_j$, the reflectability condition requires that
$\sigma_i=\sigma_j$.  Thus, for an irreducible Coxeter group, all of the
$\sigma_i$ must be equal; for a reducible Coxeter group acting on
$\R^{n_1}\oplus\R^{n_2}\oplus\cdots\oplus\R^{n_k}$, we can multiply the
coordinates in each $\R^{n_j}$ by a constant factor so that the
$\sigma_i$ are all equal.  We can then re-scale time so that all
$\sigma_i=1$.

Reflectability also requires the Coxeter group to fix the vector $\mu$
whose coordinates are the drifts $\mu_i$.  In all non-trivial cases
except for $A_{n-1}$ on $\R^n$, this requires that all $\mu_i=0$, giving
standard Brownian motion.  For $A_{n-1}$, it requires that the $\mu_i$
all be equal; we can then change coordinates to $x'_i=x_i-\mu_i t$ to
get an equivalent process in which all of the $\mu'_i$ are zero.  Thus,
if the unconstrained process has stationary increments, and is thus a
Brownian motion~\cite{harrison}, we may assume that it is standard
Brownian motion; however, we can state the theorems just as easily in
terms of the more general reflectable process.

\section{The formulas}

\begin{thm}\label{brownthm}
If $c_t$ is the density function for a reflectable continuous stochastic
process, then for absorbing boundary conditions, we have
\begin{equation}
  \btel=\sumweyl\sgn(w) c_t(w(\lambda)-\eta),\label{brownabs}
\end{equation}
and for reflecting boundary conditions, 
\begin{equation}
  \bptel=\frac{\sumweyl c_t(w(\lambda)-\eta)}
              {\#\{w\in W: w(\lambda)=\lambda\}}.\label{brownrefl}
\end{equation}
\end{thm}

If $W$ is an affine group rather than a finite group, these may be
infinite sums.  The integrals over the images of any region must
converge absolutely, since the measure of unconstrained motion over the
whole space, the set of all $W$-images of all points in the chamber, is
1.

{\em Proof.} The discrete result analogous to~\numb{brownabs} is proved
in~\cite{GZ} and~\cite{biane92}; the proof which follows is essentially
identical to that in~\cite{GZ} except that the discrete terms ``walk''
and ``step'' are replaced by their continuous analogues ``path'' and
``time.''

Every path from $\eta$ to any $w(\lambda)$ which does touch at least one
wall of the chamber has some first time $t_0$ at which it touches a
wall; let the wall be the hyperplane perpendicular to $\alpha_i$,
choosing the largest $i$ if there are several choices~\cite{Pr2}.
Reflect the path after time $t_0$ across that hyperplane; the resulting
path is a path from $\eta$ to $r_{\alpha_i} w(\lambda)$ which also first
touches wall $i$ at time $t_0$.  This clearly gives a measure-preserving
bijection of paths, and since $r_{\alpha_i}$ has sign $-1$, all such
paths cancel out in~\numb{brownabs}.  The only paths which do not cancel
in these pairs are the paths which stay within the Weyl chamber, and
since $w(\lambda)$ is inside the Weyl chamber only if $w$ is the
identity, this is the desired measure.

For~\numb{brownrefl}, we note that the map on all paths starting at
$\eta$ which takes every point to its unique image in the chamber $C$ is
measure-preserving, since we have reflecting boundary conditions and
increments which are stable under the group $W$.  This map takes all
paths which end at any $w(\lambda)$ to paths which end at $\lambda$
itself.  If $\lambda$ is on a wall of the chamber, paths to $\lambda$
may be counted multiple times, so we must divide by the size of its
stabilizer in $W$.

In practice, we can ignore this constant factor; it is 1 except for
$\lambda$ on a wall, and this is a set of measure zero unless $t=0$.
Eliminating the denominator thus changes the density function only on a
set of measure zero, and thus does not change the measure of any
measurable set.

\section{Determinant formulas for the density functions}

We can now apply this theorem to standard Brownian motion, in the Weyl
chambers of $A_{n-1}$, $B_n=C_n$, and $D_n$, with either absorbing or
reflecting boundary conditions.  The measure for unconstrained standard
Brownian motion is $c_t(\vec x)=\prodin N_t(x_i)$, where $N_t$ is the
normal distribution function with mean 0 and variance $t$.  Since this
factors into separate terms for the individual coordinates, we can use
the same techniques to compute determinant formulas as in the discrete
case~\cite{decomp}.

The most interesting case is $A_{n-1}=S_n$, the symmetric group.  The
Weyl chamber is $x_1>x_2>\cdots>x_n$.  This Brownian motion thus models
$n$ independent particles in one dimension.  With absorbing boundary
conditions, collisions are forbidden (the process terminates if one
occurs); with reflecting boundary conditions, particles collide
elastically with one another.

For absorbing boundary conditions, we write the sum~\numb{brownabs} as
\begin{equation}
  \btel=\sum_{\sigma\in S_n}\sgn(\sigma) c_t(\sigma(\lambda)-\eta)
\end{equation}
and use the value of $c_t$ to write this as
\begin{equation}
  \btel=\sum_{\sigma\in S_n}\sgn(\sigma)
          \prodin\left(N_t(\lambda_{\sigma(i)}-\eta_i)\right).
\end{equation}
This sum can be written as a determinant, which gives
\begin{equation}
  \btel=\det_{n\times n}\left|N_t(\lambda_i-\eta_j)\right|.
\label{anbrown}
\end{equation}
This determinant gives the measure for $n$ particles which start at
positions $\eta_i$ and are in independent Brownian motion to be at
positions $\lambda_i$ at time $t$ without having collided.

For a more general product of $n$ independent diffusions, with
individual density $p_t(x\to y)$ for the diffusion which started at
point $x$ to be at point $y$ at time $t$, we can apply the same argument.
This gives the following generalization of~\numb{anbrown}, which first
appears in~\cite{KM}.

\begin{equation}
  \btel=\det_{n\times n}\left|p_t(\eta_j\to\lambda_i)\right|.
\end{equation}

For reflecting boundary conditions, the calculations are the same except
that there is no term for $\sgn(\sigma)$, and thus we get permanents
rather than determinants.  This result can also be seen by observing
that an elastic collision between two identical particles is equivalent
to the two particles passing through each other with no collision, and
thus the reflected particles will be at the positions indicated by
$\lambda$ if the unreflected particles are at any permutation of the
coordinates of $\lambda$.

For $B_n$, the hyperoctahedral group, which includes permutations with
any number of sign changes, the Weyl chamber is $x_1>x_2>\cdots>x_n>0$.
This also models $n$ independent particles in one dimension, with an
additional wall at $x=0$.

We write $w\in W$ as a product of an $\epsilon$ which negates some
coordinates and a $\sigma$ in the symmetric group.  We get
\begin{equation}
  \btel(x)=\sumhypoct \sgn(\sigma)
           \prodin \epsilon_i \prodin \left(
              N_t(\epsilon_i\lambda_{\sigma(i)}-\eta_i)\right).
\end{equation}
Using the multilinearity of the products in the determinant, we can
again write this sum as a determinant, with separate terms for
$\epsilon_i=1$ and $\epsilon_i=-1$ in each entry.  We use
$N_t(x)=N_t(-x)$ to keep the signs of $\lambda$ positive and get a more
elegant formula.  This gives
\begin{equation}
  \btel=\det_{n\times n}
	  \left|N_t(\lambda_i-\eta_j)-N_t(\lambda_i+\eta_j)\right|.
\label{bnbrown}
\end{equation}
This determinant gives the measure for $n$ particles which start at
$\eta_i$ to be at $\lambda_i$ at time $t$, neither having collided nor
having touched $x=0$.  

Again, the same argument applies if $n$ particles are in general
independent diffusions, provided that the diffusions are symmetric about
$x=0$.  The more general formula is

\begin{equation}
  \btel=\det_{n\times n}
	  \left|p_t(\eta_j\to\lambda_i)-p_t(\eta_j\to-\lambda_i)\right|.
\label{bndiff}
\end{equation}

For reflecting boundary conditions, we lose the sign of the $\sigma$,
which makes the determinant into a permanent, and the sign of the
$\epsilon_i$, which turns the minus sign between the two $N_t$
in~\numb{bnbrown} or $p_t$ in~\numb{bndiff} into a plus sign.  Again,
the resulting formula is the same as would be obtained by treating
elastic collisions as though the particles passed through each other,
and allowing particles to pass through the wall at $x=0$ instead of
bouncing.  (In the transformed model, particles at positions $x$ and
$-x$ no longer collide, but since they collided elastically in the
original model and could be considered to pass through each other
instead, the effect is the same.)

For $D_n$, the even hyperoctahedral group, which includes permutations
with an even number of sign changes, the Weyl chamber is
$x_1>x_2>\cdots>x_n$, $x_{n-1}>-x_n$.  This does not give a natural
model for $n$ particles in one dimension.

Again, we can write $w=\epsilon\sigma$.  We take our sum over all
possible $\epsilon$, and then add an additional factor of
$(1+\prodin\epsilon_i)/2$ to annihilate those $\epsilon$ which are not
allowed in $D_n$.
\begin{equation}
  \btel(x)=\sumhypoct \sgn(\sigma) \frac{1+\prodin\epsilon_i}{2}
           \prodin \left(\epsilon_i 
              N_t(\epsilon_i\lambda_{\sigma(i)}-\eta_i)\right).
\end{equation}
We now take the $\frac{1}{2}$ and the $(\prodin\epsilon_i)/2$ terms
separately.  The $(\prodin\epsilon_i)/2$ term is half the sum we had
in~\numb{bnbrown}; the $\frac{1}{2}$ term gives half of~\numb{bnbrown},
but with a plus sign between the terms.  Thus we get
\begin{eqnarray}
  \btel(x) &=& \frac{1}{2}\Bigl[\det_{n\times n}
		\left|N_t(\lambda_i-\eta_j)
		-N_t(\lambda_i+\eta_j)\right|
	\nonumber\\
	   && {} + \det_{n\times n}
		\left|N_t(\lambda_i-\eta_j)
		     +N_t(\lambda_i+\eta_j)\right|\Bigr].
\label{dnbrown}
\end{eqnarray}

If we let $\lambda'$ be obtained from $\lambda$ by changing the sign of
$\lambda_n$, this will change the sign of the first term but preserve
the second term.  Thus the first term alone, with no factor of $1/2$,
is $\btel(x)-\btelp(x)$, and the second term alone is
$\btel(x)+\btelp(x)$.  If $\lambda_n=0$ or $\eta_n=0$, the
first term is zero, so $\btel(x)$ is the second term alone, with the
factor of $1/2$.

For reflecting boundary conditions, we ignore the sign of $\sigma$.
Thus the determinants become permanents, but we keep the minus signs
because they came from the factor $(1+\prodin\epsilon_i)/2$, which was
not from $\sgn(w)$.

\section{Asymptotics}

We can use these formulas to find asymptotics for the probability that
the motion will not hit a wall of the chamber by time $t$, and for its
distribution at time $t$ given that it has not hit a wall.

\subsection{Calculating the individual values for $A_{n-1}$}

We can eliminate the determinant to get a more explicit formula
for~\numb{anbrown} at a single point $\lambda$ if all the coordinates of
our starting point $\eta$ are rational.  Re-scaling by $x_i\to cx_i$,
$t\to c^2t$ will make all the coordinates integers, and we can then
translate all coordinates by $-\eta_n$ so that we have $\eta_n=0$.
(Both of these transformations leave the Weyl chamber
$x_1>x_2>\cdots>x_n$ unchanged.)

We now write out the normal distributions in~\numb{anbrown} explicitly
as exponentials, and expand $(\lambda_i-\eta_j)^2$ as
$\lambda_i^2-2\lambda_1\eta_j+\eta_j^2$: 
\begin{equation}
  \btel=\frac{1}{(2\pi t)^{n/2}}
  \det_{n\times n}\left|
       \exp\left(\frac{-\lambda_i^2+2\lambda_i\eta_j-\eta_j^2}{2t}
	   \right)\right|.
\end{equation}

Row $i$ of this matrix contains a constant factor $\exp(-\lambda_i^2)$, and
column $j$ a constant factor $\exp(-\eta_j^2)$, so we can take these
out, and put them in a constant term, which simplifies further because
$|\lambda|^2=\sum\lambda_i^2.$  This gives us
\begin{equation}
  \btel=\frac{1}{(2\pi t)^{n/2}}
	\exp\left(\frac{-|\lambda|^2-|\eta|^2}{2t}\right)
	\det_{n\times n}\left|
            \exp\left(\frac{\lambda_i\eta_j}{t}\right)\right|.
  \label{expvand}
\end{equation}

\sloppy
Since the $\eta_j$ are all integers, we can write the determinant as the
generalized Vandermonde determinant
\begin{equation}
  \det_{n\times n}\left|[\exp(\lambda_i/t)]^{\eta_j}\right|.
  \label{genvand}
\end{equation}
If $\eta_j=n-j$, this is the standard Vandermonde determinant, equal to
\begin{equation}
  \prod_{i>j}[\exp(\lambda_i/t)-\exp(\lambda_j/t)].\label{anvand}
\end{equation}
And for any non-negative integers $\eta_j$, it is the product of this
Vandermonde determinant and the Schur function~\cite{Macd}
\[
  s_{\eta_1-n+1,\eta_2-n+2,\ldots,\eta_n}
  (\exp(\lambda_1/t),\ldots,\exp(\lambda_n/t))
\]
The Schur function can also be defined combinatorially~\cite{Macd},
with the coefficient of $\prod x_i^{n_i}$ in $s_\mu$ the number of ways
to fill in the partition diagram of $\mu$, using the number $i$ exactly
$n_i$ times, such that the entries are non-decreasing in each row and
strictly increasing in each column.

\fussy
In particular, we can let $C_\eta=S_\mu(1,1,\ldots,1)$ be the total
number of such tableaux; this is important because it is the approximate
value of the Schur function when the $\lambda_i$ are much less than $t$.
(This holds because the Schur function is a homogeneous polynomial of
degree $\eta_1+\cdots+\eta_n-(n(n-1))/2$, with positive coefficients.)
This will show the dependence of the asymptotics on the starting point.

This constant is known~\cite{Macd}; it is
\begin{equation}
C_\eta=\prod_{i<j}(\eta_i-\eta_j)\bigg/\prod_{i<j}(j-i).\label{anconst}
\end{equation}
This allows us to compute the constant term in the asymptotics.

For any $\eta$ with the same sum of the coordinates, all multiples of
$1/c$, our rescaling gives a Schur function whose index is the
transformed vector $c\eta$.  Rescaling to restore the old time values
gives a Schur function which is a homogeneous polynomial in the
$\exp(\lambda_i/t)^{1/c}$ whose degree in the $\exp(\lambda_i/t)$ is
$\eta_1+\cdots+\eta_n-(n(n-1))/2c$.  Thus, for any point $\eta$, even
one with fractional coordinates, the degree will be bounded by the sum
of its coordinates.  Thus the determinants, and therefore the $\btel$,
for different starting points $\eta$ and $\eta'$ both of radius less
than a known $\delta$, will be in the approximate ratio of $C_\eta$ to
$C_{\eta'}$, with an error of $O(\delta/t)$.  If $\delta$ is fixed, then
when we calculate the asymptotics as $t\to\infty$, the Schur function
will converge to $C_\eta$ at a rate of $O(1/t)$.

\subsection{Asymptotic probability of no collisions}

The integral of the value in~\numb{expvand} over the whole Weyl chamber
is the probability that $n$ particles starting at the positions $\eta_j$
will have no collisions up to time $t$.  We can use this formula to show
that the asymptotic probability as $t\to\infty$ is a constant multiple
of $t^{-n(n-1)/4}$, with the constant depending on $\eta$.   We will also
show that, given that the $n$-dimensional Brownian motion has not hit a
wall, its radial distribution, rescaled by multiplying the radius by
$1/\sqrt{t}$, converges to the distribution of the Bessel process with
parameter $n(n+1)/2$ at time 1.

Fix $t$ very large compared with $|\eta|^2$.  If $|\lambda|^2$ is much
larger than $t$, then the exponential of $-|\lambda|^2/2t$
in~\numb{expvand} will decay exponentially fast.  In particular,
$\lambda_i/t$ can be assumed to be $O(t^{\epsilon-1/2})$, so we only
need to look at the leading nonzero terms in~\numb{genvand}; the error
will be a factor of this order when compared with the value of the
determinant.  In~\numb{anvand}, the first nonzero term of the factor
$[\exp(\lambda_i/t)-\exp(\lambda_j/t)]$ is $(\lambda_i-\lambda_j)/t$; we
can then take out the factor of $t^{-n(n-1)/2}$ and leave only a term
involving the $\lambda$.  Likewise, in the Schur function
from~\numb{genvand}, we need only keep the constant term $C_\eta$.

Thus our probability is asymptotic to the integral over the Weyl chamber of
\begin{equation}
  \frac{C_\eta}{(2\pi t)^{n/2}}
  \exp\left(\frac{-|\lambda|^2-|\eta|^2}{2t}\right)
  t^{-n(n-1)/2}\prod_{i>j}(\lambda_i-\lambda_j).
  \label{anasymp}
\end{equation}

And this integral can be computed by using Selberg's
integral~\cite{Macdroot,mehta}; we have

\begin{equation}
\int_{\R^n}\exp(-|x|^2/2) \left|\prod_{i>j}(x_i-x_j)\right|dx_i=
2^{3n/2}\prod_{k=1}^n \Gamma((k/2)+1).
\end{equation}

This corresponds to our desired integral when we set $\vec
x=\lambda/\sqrt{t}$, and divide by $n!$ because we are taking our
integral over only one of the $n!$ different Weyl chambers.  Dropping
the $\exp(-|\eta|^2/2t)$ (which goes to 1 as $t\to\infty$ and thus doesn't
affect the leading term), and writing out $C_\eta$ explicitly again gives
our full asymptotic:
\begin{equation}
\frac{\prod_{i<j}(\eta_i-\eta_j)}{\prod_{i<j}(j-i)}
(2\pi)^{n/2}
2^{3n/2}\prod_{k=1}^n[ \Gamma((k/2)+1)]
t^{-n(n-1)/4}.\label{anlead}
\end{equation}

We can also note that the exponential in~\numb{anasymp} is spherically
symmetric, while $\prod_{i<j}(\lambda_i-\lambda_j)$ is homogeneous of
degree $n(n-1)/2$.  Thus the density, integrated over the sphere of
radius $r$ at a fixed $t$, is a constant multiple of 
\[
\exp(-r^2/2t)r^{\left((n(n-1)/2\right)-1},
\]
and thus a constant multiple of the radial distribution of
$(n(n+1)/2)$-dimensional Brownian motion at time $t$.  Thus, given that
no two particles have collided, the distribution would be exactly the
same if it were given by~\numb{anasymp}.  This proves
Theorem~\ref{radialthm}; if we restrict to
$|\lambda|<ct^{1/2+\epsilon}$, the $\lambda_i/t$ terms are all
$O(t^{\epsilon-1/2})$.  Thus, in this region, the ratio of the radial
distribution for the constrained motion to the radial distribution for
unconstrained motion in $n(n+1)/2$ dimensions converges uniformly to 1
at a rate of $O(t^{\epsilon-1/2})$.

Equivalently, we could fix time and $\lambda$, and for a scalar
$\delta$, take $\hat\eta=\delta\eta$ as our starting point.  By the
scaling properties of Brownian motion, this is equivalent to keeping
$\eta$ fixed, taking $\hat\lambda=\lambda/\delta$, and $\hat
t=t/\delta^2$.  Thus, as $\delta\to 0$, the ratio of the radial
distributions converges uniformly to 1 within the region
$|\lambda|<ct^{1/2+\epsilon}$ at a rate of $O(\delta^{1-\epsilon})$, and
the probability that either distribution is outside that region goes to
0 exponentially fast.

\subsection{Asymptotics for $B_n$: no collisions and a wall}

We can use the same technique to get asymptotics for $B_n$ as for
$A_{n-1}$; here, we get a constant multiple of $t^{-n^2/2}$ as the
probability of no collision, and the radial distribution given no
collision is the $n^2+n$-dimensional Bessel process.

Here, we require that the coordinates all be odd integers.
Again, we write out the determinant~\numb{bnbrown} explicitly:
\begin{equation}
  \btel=
  \frac{1}{(2\pi t)^{n/2}}\det_{n\times n}\left|
       \exp\left(\frac{-\lambda_i^2+2\lambda_i\eta_j-\eta_j^2}{2t}
	   \right)
       -\exp\left(\frac{-\lambda_i^2-2\lambda_i\eta_j-\eta_j^2}{2t}
	   \right)\right|.
\end{equation}

As before, we remove the constant factors to get
\begin{equation}
  \btel=\frac{1}{(2\pi t)^{n/2}}
	\exp\left(\frac{-|\lambda|^2-|\eta|^2}{2t}\right)
	\det_{n\times n}\left|
        \exp\left(\frac{\lambda_i\eta_j}{t}\right)
        -\exp\left(\frac{-\lambda_i\eta_j}{t}\right)\right|.
  \label{bnexp}
\end{equation}

The determinant here is not an actual Vandermonde determinant.  However,
in the specific case $\eta_j=2n+1-2j$, we can make it a Vandermonde
determinant by elementary operations.  Adding $(-1)^k \binom{2n+1-2j}
{k}$ times column $j+k$ to column $j$ does not change the determinant,
but it changes the entries in column $j$ to
\begin{multline}
  {\sum_{k=0}^{2n+1-2j}(-1)^k \binom{2n+1-2j}{k} 
     \exp((2n+1-2j-2k)\lambda_i/t)}\nonumber \\
  =\left[\exp(\lambda_i/t)-\exp(-\lambda_i/t)\right]^{2n+1-2j}.
\end{multline}

This is a generalized Vandermonde determinant; we can make it an actual
Vandermonde determinant by dividing row $i$ by
$\exp(\lambda_i/t)-\exp(-\lambda_i/t)$.  The resulting determinant is
\begin{equation}
  \det_{n\times n}\left|[\exp(\lambda_i/t)-\exp(-\lambda_i/t)]
     ^{2(n+1-j)}\right|,
  \label{bnvand}
\end{equation}
and its value is
\begin{equation}
  \prod_{i>j}\left[(\exp(\lambda_i/t)-\exp(-\lambda_i/t))^2
             -(\exp(\lambda_j/t)-\exp(-\lambda_j/t))^2\right].
\label{bnfactor}
\end{equation}
Putting this together with the constants we have taken out, we get
\begin{eqnarray*} 
  \btel &=& \frac{1}{(2\pi t)^{n/2}}
	    \exp\left(\frac{-|\lambda|^2-|\eta|^2}{2t}\right)\\
	&&  \times
	    \left(\prod_{i=1}^n [\exp(\lambda_i/t)-\exp(-\lambda_i/t)]\right)\\
        &&  \times \prod_{i>j}\left[(\exp(\lambda_i/t)-\exp(-\lambda_i/t))^2
             -(\exp(\lambda_j/t)-\exp(-\lambda_j/t))^2\right].
\end{eqnarray*}

For a more general starting point, with all coordinates odd integers, we
note that $(\exp(\lambda_i\eta_j/t)-\exp(-\lambda_i\eta_j/t))$ is a
polynomial in $\exp(\lambda_i/t)-\exp(-\lambda_i/t)$ with no constant
term.  We can break these polynomials into their individual terms,
giving a large number of determinants, each one a generalized
Vandermonde determinant in the $\exp(\lambda_i/t)-\exp(-\lambda_i/t)$.
Each individual determinant is thus the product of~\numb{bnfactor} and a
Schur function; it also contains the product of the
$\exp(\lambda_i/t)-\exp(-\lambda_i/t)$ as a factor, since these are
constant factors in row $i$.  Thus each determinant is a product of
these factors with some symmetric function in the 
$\exp(\lambda_i/t)-\exp(-\lambda_i/t)$.  As with $A_{n-1}$, the error we
get in approximating the Schur functions by their constant term $C_\eta$
is a factor of $O(1/t)$.

And as with $A_{n-1}$, we can get asymptotics by integrating this over the
Weyl chamber.  For large $t$, the radial exponential will be
exponentially small if $|\lambda|>t^{1/2+\epsilon}$, so we can assume
that all of the $\lambda_i/t$ are very small.  Thus
$\exp(\lambda_i/t)-\exp(-\lambda_i/t)$ can be approximated by its
leading nonzero term, $2\lambda_i/t$.  Thus, as with~\numb{anasymp}, we
get
\begin{equation}
  \frac{1}{(2\pi t)^{n/2}}   
  \exp\left(\frac{-|\lambda|^2-|\eta|^2}{2t}\right)
  (2/t)^{n^2}\prodin(\lambda_i)\prod_{i>j}(\lambda_i^2-\lambda_j^2).
  \label{bnasymp}
\end{equation} 

This integral can also be computed by using Selberg's
integral~\cite{Macdroot,mehta}; we have
\begin{multline*}
\lefteqn{\int_{\R^n}\exp(-|x|^2/2) \prodin|x_i|
\left|\prod_{i>j}(x_i^2-x_j^2)\right| dx_i}\\
=\frac{2^{(n^2+3n)/2}}{\pi^{n/2}}
\prodin\left[\Gamma(1+i/2)\Gamma((1+i)/2)\right].
\end{multline*}

This corresponds to our desired integral when we set $\vec
x=\lambda/\sqrt{t}$, and divide by $2^n n!$ because that is the number
of Weyl chambers.  Dropping the $\exp(-|\eta|^2/2t)$ because it goes to 1
as $1/t$, we get our asymptotic probability that Brownian motion started
at the specific point $\eta$ will not hit a wall up to time $t$:
\begin{equation}
  \frac{2^{(3n^2/2)}}{\pi^n n!}
  \prodkn\left[\Gamma(1+k/2)\Gamma((1+k)/2)\right]
  t^{-n^2/2}.\label{bnlead}
\end{equation}

For a more general starting point, each individual determinant in the
sum gives an integral of a homogeneous polynomial in the
$\exp(\lambda_i/t)-\exp(-\lambda_i/t)$ which is of degree at least
$n^2$, since it contains the previous determinant as a factor.  The same
technique as above gives an asymptotic which is thus at most
$t^{-n^2/2}$.  Thus the asymptotic is at most a constant multiple of
$t^{-n^2/2}$, and is less if and only if the coefficient of $t^{-n^2/2}$
is zero.  But the coefficient cannot be zero.  Let
$m=\max((\eta_i-\eta_{i+1}/2),\eta_n)$, and $\eta'_i=(2n+1-2i)m$.  Then
we know that Brownian motion starting at $\eta'$ decays asymptotically
as $t^{-n^2/2}$, but Brownian motion starting at $\eta$ will always hit
a wall if it is translated to start at $\eta'$, because $\eta'$ is at
least as far as $\eta$ from every wall.  Thus the asymptotic probability
will be some constant $C_\eta$ times the formula~\numb{bnlead}; we will
compute $C_\eta$ in Section~\ref{constsec}.

A result analogous to Theorem~\ref{radialthm} also holds for $B_n$.
The exponential in~numb{bnasymp} is spherically symmetric, while the
product is homogeneous of degree $n^2$.  The same argument as for
$A_{n-1}$ shows that the ratio between the radial part of this
distribution and the radial distribution for unconstrained Brownian
motion in $n^2+n$ dimensions converges uniformly to 1 at a rate of
$O(t^{\epsilon-1/2})$ inside the region $|\lambda|<ct^{1/2+\epsilon}$,
and the probability that either distribution is outside that region goes
to 0 exponentially fast.  Again, we can rescale by multiplying the
radius by $1/\sqrt{t}$ to get convergence to a fixed distribution.

\subsection{Asymptotics for $D_n$}
The process for $D_n$ is almost the same as for $B_n$, so we won't work
it out in full detail; we get a similar result with the same error
terms and convergence properties.

Here, it is most natural to let $\eta_i=n-i$.  For this value of $\eta$,
we have only the second determinant in~\numb{dnbrown}, with a plus sign
between the terms; the last row of the other determinant is zero.  For
general $\eta$, the first determinant is the $B_n$ determinant, which we know
is $O(t^{-n^2/2})$, and we will show that the second determinant is
asymptotically larger.

The determinant that we get is
\begin{equation}
  \det_{n\times n}\left|
        \exp\left(\frac{\lambda_i\eta_j}{t}\right)
        +\exp\left(\frac{-\lambda_i\eta_j}{t}\right)\right|.
  \label{dnexp}
\end{equation}

Again, this isn't a Vandermonde determinant, but
$\exp(\lambda_i/t)+\exp(-\lambda_i/t)^n$ is equal to
$\exp(\lambda_i\eta_j/t)+\exp(-\lambda_i\eta_j/t)$ plus
a sum of lower order terms, so elementary operations which do not
change the determinant give us the Vandermonde determinant
\begin{equation}
  \det_{n\times n}\left|[\exp(\lambda_i/t)+\exp(-\lambda_i/t)]^{n-j}\right|,
  \label{dnvand}
\end{equation}
and its value is
\begin{equation}
  \prod_{i>j}\left[\exp(\lambda_i/t)+\exp(-\lambda_i/t)
             -\exp(\lambda_j/t)-\exp(-\lambda_j/t)\right].
\label{dnfactor}
\end{equation}

The leading nonzero term is
\[
  \prod_{i>j}\frac{(\lambda_i/t)^2-(\lambda_j/t)^2}{2},
\]
and the full integral is

\begin{equation}
  \frac{1}{(2\pi t)^{n/2}}   
  \exp\left(\frac{-|\lambda|^2-|\eta|^2}{2t}\right)
  (2/t)^{n^2-n}\prod_{i>j}(\lambda_i^2-\lambda_j^2).
  \label{dnasymp}
\end{equation} 

Again, we get a result which can be obtained from Selberg's
integral~\cite{Macdroot,mehta}; we have
\begin{equation}
\int_{\R^n}\exp(-|x|^2/2)
\left|\prod_{i>j}(x_i^2-x_j^2)\right| dx_i=
\frac{2^{(n^2+2n)/2}}{\pi^{n/2}}
\prodin\left[\Gamma(1+i/2)\Gamma(i/2)\right].
\end{equation}

Here, there are $2^{n-1}n!$ Weyl chambers, which gives the asymptotic
\begin{equation}
  \frac{2^{(3n^2-3n+2)/2}}{\pi^n n!}
  \prodkn\left[\Gamma(1+k/2)\Gamma((1+k)/2)\right]
  t^{(-n^2+n)/2}.\label{dnlead}
\end{equation}

For a general starting point whose coordinates are all integers, we use
the same technique as for $B_n$.  The terms in the determinant are all
polynomials in $\exp(\lambda_i/t)+\exp(-\lambda_i/t)$, so we can again
split the sum into individual determinants, each of which is the product
of~\numb{dnfactor} and a Schur function of the
$\exp(\lambda_i/t)+\exp(-\lambda_i/t)$.  Since $\lambda_i/t$ is small,
these are all close to 2, so the full asymptotic is a sum of terms of
order $t^{(-n^2+n)/2}$.  As with $B_n$, we can translate to $\eta'$ which
is further from any wall than $\eta$ to show that the coefficient of
$t^{(-n^2+n)/2}$ cannot be zero.  

\subsection{The constant factor for a general starting
point}\label{constsec}

Since we know the constant factor in the asymptotic probability of no
collision up to time $t$ for one specific starting point $\eta$, and the
asymptotic distribution for an arbitrary starting point, we can use the
time-reversibility of Brownian motion to compute the asymptotic
probability of no collisions for an arbitrary starting point.  The
argument is the same for all of the Weyl groups.

In each case, the density $b_t(\eta,\lambda)$ for Brownian motion
starting at $\eta$ to be at $\lambda$ at time $t$, not having collided
with a wall up to time $t$, is asymptotic for large $t$ to a product of 
the form 
\begin{equation}
C_\eta C'_\lambda f(t) 
  \exp\left(\frac{-|\lambda|^2-|\eta|^2}{2t}\right), \label{genasymp}
\end{equation}
in which $C'_\lambda$ and $f(t)$ are known.  Since Brownian motion is
symmetric in time, the density $b_t(\lambda,\eta)$ must be equal to
$b_t(\eta,\lambda)$.  As long as $t$ is large enough compared to
$|\eta|$ and $|\lambda|$ for the formula~\numb{genasymp} to be valid
(which it will be for large $t$ because of the exponentials), we can
reverse the roles of $\eta$ and $\lambda$.  Thus we have $C_\eta
C'_\lambda=C_\lambda C'_\eta$, and since the formula for $C'$ is known,
we see that $C_\eta$ must be a constant multiple of $C'_\eta$; it is
thus a constant multiple of $\prod (\eta_i-\eta_j)$ for $A_n$, of $\prod
(\eta_i^2-\eta_j^2)\prod\eta_i$ for $B_n$, and $\prod
(\eta_i^2-\eta_j^2)$ for $D_n$.

We know the value of the constant from the formulas \numb{anlead},
\numb{bnlead}, and \numb{dnlead}.  For $A_n$, we already have the
value for a general starting point because we used the Schur functions
to compute $C_\eta$ in~\numb{anlead}; we could have instead used this
technique. For $B_n$, we get
\begin{eqnarray*}
  \prod_{i<j}\frac{\eta_i^2-\eta_j^2}{(2j-1)^2-(2i-1)^2}
  \prod_{i=1}^n\frac{\eta_i}{2n+1-2i}\nonumber\\
  \quad\times\frac{2^{(3n^2/2)}}{\pi^n n!}
  \prodkn\left[\Gamma(1+k/2)\Gamma((1+k)/2)\right] t^{-n^2/2}.
\end{eqnarray*}
as the asymptotic probability of no collision; for $D_n$, we get
\begin{eqnarray*}
  \prod_{i<j}\frac{\eta_i^2-\eta_j^2}{(2j-1)^2-(2i-1)^2}\nonumber\\
  \quad\times  \frac{2^{(3n^2-3n+2)/2}}{\pi^n t^{n^2/2}n!}
  \prodkn\left[\Gamma(1+k/2)\Gamma((1+k)/2)\right] t^{(-n^2+n)/2}.
\end{eqnarray*}

\section{Random matrices and conditioned Brownian motion}

\subsection{Brownian motion on a Lie algebra}

Instead of viewing our Weyl chamber as the chamber of a Weyl group, we
can view it as the space of eigenvalues of the Lie algebra corresponding
to that Weyl group, and then use the theory of Lie algebras to study it.
The following construction was first developed by Dyson~\cite{dyson} for
$\su_n$; he computed the properties of the Brownian motion from the
specific data rather than using the general Lie theory.  The results we
use from Lie theory are given in~\cite{adams,Hu}.

Given a Lie group $G$ and its Lie algebra $\mathfrak g$, we can define a
normal distribution or Brownian motion on $\mathfrak g$ by viewing it as a
vector space and taking an orthonormal basis with respect to a
$G$-invariant inner product.  In particular, on $\so_n(\R)$, the Lie
algebra of skew-symmetric matrices, we can take $1/\sqrt{2}$ times standard
Brownian motion on each matrix entry $M_{ij}$ with $i<j$, and take
$M_{ji}=-M_{ij}$; this makes the Hilbert-Schmidt norm of $M$ equivalent
to the radius of a Brownian motion in $n^2-n$ dimensions.  Likewise, for
$\su_n(\C)$, the Lie algebra of skew-Hermitian matrices, we can take
standard Brownian motion on the imaginary part of each diagonal entry
$M_{ii}$, and $1/\sqrt{2}$ times standard Brownian motion on each matrix
entry $M_{ij}$ with $i<j$, with $M_{ji}=-\overline{M_{ij}}$; this again
gives us the Hilbert-Schmidt norm, this time in $n^2$ dimensions.  For
any Lie algebra, the Brownian motion at time $t$ will have a Gaussian
distribution on the space of matrices in that Lie algebra.  This allows
us to use all of the known results about random matrices~\cite{mehta} to
study the motion, and its eigenvalues in particular.

Since this Brownian motion is invariant under conjugation by $G$, it
induces a diffusion on the space of eigenvalues of the matrices.  The
Weyl chamber is in a correspondence with the space of eigenvalues,
obtained by dividing the independent eigenvalues by $i$ and arranging
them in decreasing order.  To study this diffusion, we consider the
measure on the Weyl chamber which is induced by the standard measure on
the Lie algebra.  We can study this measure by proving a version of the
Weyl Integration Formula for Lie algebras, which can be proved directly
from the analogous formula for Lie groups~\cite{adams}.

\begin{lemma}
Let $\mathfrak g$ be a finite-dimensional Lie algebra with $m$ roots, and
$f(X)=\hat f(\lambda_1,\ldots,\lambda_n)$ a function on $\mathfrak g$ with
compact support which depends only on the eigenvalues of $X$.  Let
$\delta'$ be $(2\pi)^m$ times the product of all the roots of $\mathfrak g$.
Then the integral of $f$ over $\mathfrak g$ is equal to the integral of 
$\delta'\bar\delta'\hat f$ over the Weyl chamber.
\end{lemma}

{\em Proof.} Let $G$ be the Lie group corresponding to $\mathfrak g$.  We
can apply the exponential map $h: X\to \exp(\pi itX)$ to translate a
neighborhood of zero in $\mathfrak g$ which contains the support of $f$ to
an arbitrarily small neighborhood of the origin in $G$ by choosing $t$
sufficiently small.  The derivative of this map can be made arbitrarily
close to $(\pi it)^n$ (which becomes $(\pi t)^n$ when we consider volume
elements) by taking $t$ sufficiently small.  We now apply the Weyl
Integration Formula to $h\circ f$, which is a class function on the Lie
group $G$; it says that the integral of $h\circ f$ over $G$ is equal to
the integral of $\delta\bar\delta$ over a maximal torus, where 
\begin{equation}
\delta=\prod_{\alpha\in\Delta}\left(e^{\pi i \alpha(\vec x)}
-e^{-\pi i\alpha(\vec x)}\right).
\end{equation}
For $\vec x$ sufficiently close to zero, we can approximate $\delta$ by
its leading nonzero term, which is 
\begin{equation}
\delta\approx\delta'=(2\pi)^m\prod_{\alpha\in\Delta}\alpha(\vec x).
\end{equation}
This allows us to proceed with the integration.

We can apply this formula to the characteristic functions of
small intervals for the eigenvalues; this shows that the measure at a
point $(x_1,\ldots,x_n)$ in the Weyl chamber is proportional to
$\delta\bar\delta$.  In particular, we note that $\delta$ is positive real
on the interior of the Weyl chamber, and zero on the walls; thus we have
a diffusion which is conditioned never to leave the Weyl chamber.

Note that this proves Theorem~\ref{processthm} for this process, and an
analogous theorem for any Weyl group.  If we view $\mathfrak g$ as a vector
space, the Hilbert-Schmidt norm of a matrix is the radius of a vector.
If there are $m$ positive roots, then $\mathfrak g$ has dimension
$n+2m$~\cite{Hu}.  The Hilbert-Schmidt norm is thus given by the Bessel
process in $n+2m$ dimensions.  The Hilbert-Schmidt norm is unchanged
under diagonalization, so the diffusion in the Weyl chamber has its
radius given by the same process.

To properly state this result for Brownian motion conditioned never to
leave the Weyl chamber, we cannot allow the process to start at the
origin, which is not a point in the chamber.  However, the square of the
product of the roots is a homogeneous function of degree $2m$; it thus
gives an identical distribution on any fixed radius.  Thus, if we start
the motion on the Lie algebra at 0, then at any later time, given the
fixed radius, the distribution on the sphere will be given by the square
of the product of the roots, and the process can be continued from that
time on as a Bessel process.  This allows us to state the general
theorem, formalizing Theorem~\ref{processthm}.

\begin{thm}
For any Weyl group $W$ acting on $\R^n$ with $m$ roots, consider the
process which starts at a fixed radius $r_0$, with the starting point
chosen on the sphere of radius $r_0$ by a distribution with density
proportional to the square of the product of the roots.  Then the radius
of the position of transformed motion at time $t$ gives a Bessel process
with parameter $n+2m$.
\end{thm}

The simplest case of this theorem is for $B_1$, which is one-dimensional
motion with an absorbing boundary at $x=0$.  In this case, the only root
is $x$, so $m=1$, and we have the result of~\cite{pitman}
and~\cite{williams} that Brownian motion conditioned never to hit 0 is
the same as the three-dimensional Bessel process.  In our more general
cases, the number of dimensions for the Bessel process is $n+2m=n^2$
for~$A_{n-1}$, $n(2n+1)$ for~$B_n$, and $n(2n-1)$ for~$D_n$.

\subsection{Construction by $h$-transformation}

We could have alternatively constructed this diffusion by
$h$-transformation~\cite{durett}.  We will first construct the
conditioned Brownian motion in this way, and then show that the two
processes are actually identical; this allows us to develop the physical
model of Theorem~\ref{forcethm}.

\subsubsection{General properties}

Given standard Brownian motion in any number of dimensions, we can use
the process of $h$-transformation to construct a Brownian motion
satisfying certain conditions.  For any non-negative harmonic function
$h$, the measure for the transformed Brownian motion to go to $\lambda$
at time $t$ after starting at $\eta$ is $h(\lambda)$ times the measure
for untransformed Brownian motion with the same starting point,
renormalized by an appropriate constant so that the total measure on all
paths is 1.

In particular, suppose that we have a connected region $D$ and a
harmonic function $h$ which is zero on the boundary of the region and
positive on the interior.  Then the $h$-transformed Brownian motion will
be conditioned to stay in the interior of the region.

The density function for this transformed motion will be $h$ times the
density function for the untransformed motion on the same region,
normalized appropriately.  That is, its value at a point $\lambda$ and
time $t$ will be $h(\lambda)$ times the value $\btel$ which gives the
measure for unconstrained motion to go from $\eta$ to $\lambda$ in time
$t$ while staying within the region, since the transformation puts a new
measure on the same set of paths.

The same technique can also be applied to a discrete random walk with
a set $S$ of steps, provided that $h$ is harmonic in the discrete
lattice; that is, we need
\begin{equation}
\sum_{s\in S} h(x+s)/|S|=h(x).
\end{equation}
The transformed discrete walk now has probability $h(x+s)/(h(x)|S|)$
instead of $1/|S|$ of going from $x$ to $x+s$ in a given step.  Thus,
given a starting point $\eta$, the probability of going to $\lambda$ in
a given number of steps of the transformed random walk is
$h(\lambda)/h(\eta)$ times the probability for the untransformed walk to
go to $\lambda$ while staying within the region in which $h$ is
positive.  

In order to make it impossible for the walk to leave the region $D$, we
need $h$ to be zero on all points which can be reached from the interior
of $D$ in a single step.  Thus, only if the discrete walk is reflectable
(as defined above) is it sufficient for $h$ to be zero on the continuous
boundary of $D$.  

\subsubsection{Finding the function}

The properties of conditioning, as well as our asymptotics, suggest that
$\delta$ itself should be our $h$.  It can be checked algebraically that
$h$ is harmonic for each group; however, it can also be proved naturally.

\begin{thm}
For any finite Coxeter group $W$, the product of all the positive roots
is a harmonic function, both for the continuous Laplacian and for the
discrete Laplacian
\[
L_S h(\vec x) = \frac{1}{|S|}\sum_{\vec s\in S}
                   [h(\vec x+\vec s)-h(\vec x)],
\]
for any set $S$ which is symmetric under the group $W$.
\end{thm}

{\em Proof.}  This result in the continuous case is due to~\cite{BB};
the discrete argument is a simple generalization which is mentioned
in~\cite{biane95}.

By the properties of root systems, a reflection in any
simple root changes the sign of only that root, while permuting the
other positive roots.  Thus the product $h$ of all $m$ roots is
antisymmetric in all the simple roots, and by applying repeated
reflections, we see that it is antisymmetric in every root.  It is also
of degree $m$.  

The continuous Laplacian is spherically symmetric, and thus symmetric
under $W$.  The discrete Laplacian is symmetric under $W$ because the
set $S$ is.  Applying the continuous Laplacian to a polynomial decreases
the degree by 2, while applying the discrete Laplacian decreases the
degree by at least 1 since $h(\vec x+\vec s)-h(\vec x)$ is of lower
degree than $h$.  Thus the application of the Laplacian to $h$ gives a
polynomial which is of degree less than $m$ which is still antisymmetric
in $W$.

Now, any polynomial which is antisymmetric in $W$ must be zero on every
one of the hyperplanes of reflection.  If it is not identically zero, it
must have all $m$ roots as factors, so it must be of degree at least
$m$.  Thus the Laplacian must annihilate our polynomial $h$, so $h$ is
harmonic for either the discrete or continuous walk.

It follows that the transformation by this $h$ gives the same process as
the process generated by Lie theory.  We have already computed the
asymptotic density that Brownian motion started at a fixed $\eta$ will
remain in the Weyl chamber for time $t$ and be at $\lambda$ at that
time.  Transforming by $h$ has the effect of multiplying the measure of
all paths from a fixed $\eta$ to an arbitrary $\lambda$ which stay
within the chamber by a factor of $h(\lambda)$ (and a normalizing
constant).  Thus the density function for the transformed Brownian
motion which starts at $\eta$ to be at $\lambda$ at time $t$ is
proportional to $h(\lambda)\btel$, and this converges to a constant
multiple of $h(\lambda)^2$ as $t$ becomes large.  This is the same
factor $\delta'\bar\delta'$ which we obtained from the Weyl Integration
Formula.

\subsection{The infinitesimal generator and physical models}

We have constructed the conditioned process as a transformation of
Brownian motion.  We can also construct it as a diffusion with its known
infinitesimal generator.  In this form, both the multidimensional models
and the models of $n$ independent particles lead to natural physical
models. 

We will use the notation of~\cite{KT} for infinitesimal generators of
diffusions.  The drift vector $\mu$ is defined by
\[
\mu_i(\vec x, t) = \lim_{\Delta t\downarrow 0} \frac{1}{\Delta t}
  \{ E(X_i(t+\Delta t))-X_i(t) | \vec X(t)=\vec x\},
\]
and the infinitesimal variance matrix is defined by 
\[
\sigma_{i,j}(\vec x, t) = \lim_{\Delta t\downarrow 0} \frac{1}{\Delta t}
  \{ [E(X_i(t+\Delta t))-X_i(t)][E(X_j(t+\Delta t))-X_j(t)] | \vec
  X(t)=\vec x\},
\]
We will omit the variable $t$ in the infinitesimal generators of our
diffusions, because they are independent of time.

It can easily be checked that Brownian motion transformed by a harmonic
(and thus necessarily $C^2$) function $h$ has infinitesimal drift
$\mu=\nabla h/h$ and infinitesimal variance $\sigma_{ij}=\delta_{ij}$.

Note that $\nabla h/h=\nabla(\log h)$.  This makes the computation
easy, because our $h$ is the product $\prod_{\alpha\in\Delta}\alpha$ of
all the roots, viewed as linear functions of the $x_i$.  If we write
$\alpha(\vec x)$ as the dot-product $(\alpha,\vec x)$. then we have
\begin{equation}
\frac{\nabla h}{h} = \sum_{\alpha\in\Delta}\nabla(\alpha,\vec x)
= \sum_{\alpha\in\Delta}\alpha\frac{1}{(\alpha,\vec x)}.
\end{equation}

Since $(\alpha,\vec x)/|\alpha|$ is the distance from $\vec x$ to the
hyperplane orthogonal to $\alpha$, while $\alpha/|\alpha|$ is the unit
vector in the direction of $\alpha$, this term in the drift is the
inverse of the distance between $\vec x$ and the hyperplane, directed
away from the hyperplane.  We thus have the following physical model.

\begin{thm}\label{repelthm}
For any finite Coxeter group, Brownian motion conditioned to stay within
a chamber is equivalent to the motion of a particle in a diffusion with
constant infinitesimal velocity, and a repulsive force from every
hyperplane of reflection (not merely the walls of the chamber) inversely
proportional to the distance from that hyperplane.
\end{thm}

We can also look at the $n$ coordinates as individual motions in one
dimension.  If the root $\alpha$ contains $cx_i$, we get a term
$c/(\alpha,\vec x)$ in the sum.

In particular, for $A_n$, the drift $\mu_i$ is $\sum_{i\ne
j}1/(x_i-x_j)$.  Thus each particle is subject to a repelling force from
every other particle (not merely its neighbors), inversely proportional
to the distance between them.  This proves Theorem~\ref{forcethm}; this
result is originally due to Dyson~\cite{dyson}.  

The model for $B_n$ is not as natural as a model of particles.  The
drift $\mu_i$ is 
\[
\frac{1}{x_i}
+\sum_{i\ne j}\left(\frac{1}{x_i-x_j}+\frac{1}{x_i+x_j}\right)
\]
That is, each particle is repelled by every other particle, and by the
wall at 0 (the $1/x_i$ term), but also by the mirror image of every
other particle (the $1/(x_i+x_j)$ term), as if the wall at $x=0$ was
also a mirror reflecting all forces.  For $D_n$, we have only the terms
of $1/(x_i-x_j)$ and $1/(x_i+x_j)$; this means that the mirror reflects
forces but is itself permeable to particles.  

These models are more natural as models of the eigenvalues.  For $D_n$,
a matrix in the Lie algebra $\so_{2n}$ has eigenvalues $\pm i\lambda_j$;
this model thus says that the eigenvalues in different pairs $\pm
i\lambda_j$ and $\pm i\lambda_k$ for $j\ne l$ repel one another other,
although the pair $\pm i\lambda_j$ do not repel each other.  For $B_n$,
a matrix in the Lie algebra $\so_{2n+1}$ has eigenvalues $\pm
i\lambda_j$ and $0$; this model thus says that the unconnected
eigenvalues $i\lambda_j$ and $i\lambda_k$ for $j\ne k$ repel each other,
and each eigenvalue is also repelled by the fixed eigenvalue at 0.

\section{Open problems}

The discussion of Brownian motion on a Lie algebra is valid for the
exceptional Lie groups as well, but the techniques for computing the
specific asymptotics do not appear to work.  For a chamber of a general
Coxeter group, the Lie algebra technique is not meaningful.  In either
case, is it possible to get the same type of asymptotics, with constant
terms in particular?

We have shown that the same harmonic functions which we used to
transform Brownian motion can be used to transform discrete random
walks.  Is it possible to use these results to compute asymptotics for
the discrete walks, including the constants on leading terms?  

Brownian motion can be defined on a general manifold~\cite{stochast}.
This allows us to apply the argument of Theorem~\ref{brownthm} whenever
we have a suitable chamber.  As before, the Brownian motion must be
symmetric under the reflections in any wall of the chamber, and the
reflections in the walls must generate a discrete group which partitions
the manifold into chambers.  For example, since 2-dimensional Brownian
motion is conformally invariant, we can define a Brownian motion on the
modular surface~\cite{stochast}.  Our chamber can be the standard
fundamental domain; if we use the standard map of the modular surface to
the upper half-plane, our chamber is bounded by $x>-1/2$, $x<1/2$, and
$x^2+y^2>1$.  Can the resulting formulas be used to compute properties
of this Brownian motion, such as asymptotic survival probabilities,
hitting times, and physical models?

\medbreak 

{\bf Acknowledgements.}  Sections 2-4 are part of my thesis, ``Walks and
Representation Theory,'' which was written at Harvard University.  I
would like to thank Peter Magyar for his collaboration in work on the
discrete problem; Itai Benjamini and Gil Kalai, for suggestions related
to the asymptotic results; Arun Ram, for a suggestion which simplified
the computation of the determinant formulas; Dan Stroock, for suggesting
that I study reflecting boundary conditions; the referee of an earlier
version of this paper, for suggesting the Lie theory argument and
connection to random matrices; Wendelin Werner, for suggesting some
important references; my thesis advisors, Persi Diaconis and Richard
Stanley, for their guidance in my thesis research; and Yuval Peres, for
his guidance in the research beyond my thesis.


\begin{thebibliography}{99}

% \parskip=\baselineskip
% \ifjournal\leavevmode\vskip-\parskip\fi % To get indentation right

\bibitem {adams} Adams, J. F., {\em Lectures on Lie
Groups.}  University of Chicago Press, 1969.

%bibitem S. Axler, P. Bourdon, and W. Ramey, {\em Harmonic
%Function Theory.} Springer-Verlag, New York, 1992.

\bibitem{BB} B\'erard, P., and Besson,
G., Spectres et groupes cristallographiques II: domaines
sph\'eriques. {\em Ann. Inst. Fourier} {\bf 30}(1980):3, 237--248.

\bibitem{biane92} Biane, P., Minuscule weights and
random walks on lattices. {\em Quant. Prob. Rel. Topics} {\bf 7}(1992)
51--65. 

\bibitem{biane} Biane, P., Quelques propriet\'es du
mouvement Brownien dans un cone. {\em Stoch. Proc. Appl.} {\bf
53}(1994) 233--240.

\bibitem{biane95} Biane, P., Permutation model for
semi-circular systems and quantum random walks. {\em Pacific J. Math.}
{\bf 171}(1995) 373--387.

\bibitem{Bou} Bourbaki, N., {\em Groupes et Algebres
de Lie,} Chapters 4,5,6. Hermann, Paris, 1968.

\bibitem{burkholder} Burkholder, D. L., Exit times
of Brownian motion, harmonic majorization, and Hardy spaces. {\em
Adv. Math.} {\bf 26}(1977) 182--205.

\bibitem{deblassie} DeBlassie, R. D.,  Exit times
from cones in $\R^n$ of Brownian motion. {\em Prob. Th. Rel. Fields\/}
{\bf 74}(1987) 1--29.

\bibitem{durett} Durett, R. {\em Brownian Motion and
Martingales in Analysis.} Wadsworth, Belmont, California, 1984.

\bibitem{dyson} Dyson, F. J.. A Brownian-motion model
for the eigenvalues of a random matrix, {\em J. Math. Phys.}{\bf
3}(1962), 1191--1198.

\bibitem{freedman} Freedman, D., {\em Brownian Motion and Diffusion.}
Holden-Day, San Francisco; second edition published by Springer-Verlag,
New York, 1971)

\bibitem{GZ} Gessel, I. M., and Zeilberger, D.  Random walk in a Weyl
chamber. {\em Proc. Amer. Math. Soc.}  {\bf 115}(1992), 27--31.

%bibitem D. J. Grabiner, ``A combinatorial correspondence
%for walks in Weyl chambers,'' {\em J. Combin. Th. A}, {\bf 71}(1995),
%275--292.

\bibitem{decomp} Grabiner, D. J., and Magyar,
P.  Random walks in Weyl chambers and the decomposition of tensor
powers. {\em J. Alg. Combin.}  {\bf 2}(1993) 239--260.

\bibitem{harrison} Harrison, J. M. {\em Brownian
Motion and Stochastic Flow Systems.} John Wiley and Sons, New York, 1985.

\bibitem{Hu} Humphreys, J. E. {\em Introduction to
Lie Algebras and Representation Theory.} Springer-Verlag, New York, 1972.

\bibitem{HW} Hobson, D. and Werner, W. Non-colliding Brownian motion on
the circle. {\em Bull. London. Math. Soc.}{\bf 28}(1996), 643--650.

\bibitem{KM} Karlin, S. P., and MacGregor,
G. Coincidence probabilities. {\em Pacific. J. Math.} {\bf
9}(1959),1141--1164.

\bibitem{KT} Karlin, S. P., and Taylor,
H. M., {\em A Second Course in Stochastic Processes.} Academic
Press, New York, 1981.

\bibitem{Macd} Macdonald, I. G. {\em Symmetric
Functions and Hall Polynomials.} Oxford University Press, 1979.

\bibitem{Macdroot} Macdonald, I. G.  Some conjectures for root
systems. {\em SIAM J. Math. Anal.} {\bf 13}(1982), 988-1007.

\bibitem{stochast} McKean, H. P.  {\em Stochastic
Integrals.} Academic Press, New York, 1969.

\bibitem{mehta} Mehta, M. L.  {\em Random Matrices:
Second Edition}.  Academic Press, New York, 1991.

%bibitem T. Muir (revised by W. H. Metzler), {\em A Treatise on
%the Theory of Determinants}, Dover Publications, New York, 1928.

\bibitem{OU} O'Connell, N. and Unwin,
A.  Collision times and exit times from cones: A duality. {\em
Stoch. Proc. Appl.} {\bf 43}(1992), 291--301.

%bibitem R. G. Pinsky, ``On the convergence of diffusion
%processes conditioned to remain in a bounded region for large times to
%limiting positive recurrent diffusion processes,'' {\em Ann. Prob.} {\bf
%13}(1985), 363--378.

\bibitem{pitman} Pitman, J. One-dimensional Brownian
motion and the three-dimensional Bessel process. {\em Ann. Appl. Prob.}
{\bf 7}(1975), 511--526.

%bibitem G. Polya, ``Qualitatives \"uber Warmeausgliechung,''
%{\em Z. Angew. Math. Mech.} {\bf 13}(1933), 125--128.

% This citation used for consistency with thesis
\bibitem{Pr2} Proctor, R. A.  Reflection and algorithm
proofs of some more Lie group dual pair identities. {\em
J. Combin. Th. A\/} {\bf 62}(1993) 107--127.

\bibitem{selberg} Selberg, A. Bemerkinger om et
multipelt integral. {\em Norsk Mathematisk Tidsskrift\/} {\bf 26}(1944)
71-78.

\bibitem{WM} Watanabe, T. and Mohanty, S. G., On an inclusion-exclusion
formula based on the reflection principle. {\em Discrete Math} {\bf
64}(1987), 281--288.

%bibitem H. Wilf, ``Ascending subsequences of permutations and
%the shapes of tableaux,'' {\em J. Combin. Th. A\/} {\bf 60}(1992),
%155--157.

\bibitem{williams} Williams, D. (1974). Path decomposition
and continuity of local time for one-dimensional diffusions. {\em
Proc. London Math. Soc.} {\bf 28}(1974), 738--768.

\bibitem{Z} Zeilberger, D. Andre's reflection
proof generalized to the many-candidate ballot problem. {\em Discrete
Math} {\bf 44}(1983) 325--326.

\end{thebibliography}
\end{document}